\documentclass[11pt,english]{article}
\usepackage[T1]{fontenc}
\usepackage[latin9]{inputenc}
\usepackage[a4paper]{geometry}
\usepackage{array}
\usepackage{float}
\usepackage{multirow}
\usepackage{amsmath}
\usepackage{amssymb}
\usepackage{graphicx}

\makeatletter

\providecommand{\tabularnewline}{\\}

\usepackage{amsmath}
\usepackage{times}
\numberwithin{equation}{section}
\newtheorem{The}{Theorem}[section]
\newtheorem{Lem}{Lemma}[section]

\newtheorem{Exp}{Example}[section]

\def\bproof{\textbf{Proof}: }

\def\eproof{\hfill$\Box$}

\date{}
\let\oldnormalsize=\normalsize
\renewcommand{\normalsize}{%
  \oldnormalsize%
  \setlength{\abovedisplayskip}{5pt}%
  \setlength{\abovedisplayshortskip}{5pt}%
  \setlength{\belowdisplayskip}{5pt}%
  \setlength{\belowdisplayshortskip}{5pt}}
\normalsize
\usepackage{capt-of}

\makeatother

\usepackage{babel}
\begin{document}
\title{Higher order graded mesh scheme for time fractional differential equations}
\author{G. Naga Raju$^{*}${\normalsize{}$\;,$} H. Madduri$^{*}$\footnote{gnagaraju@mth.vnit.ac.in, **harshi.madduri@gmail.com}\\
{\footnotesize{}Department of Mathematics, Visvesvaraya National Institute
of Technology, Nagpur, India.}}

\maketitle
\vspace{-0.5cm}

\begin{abstract}
In this article, we propose a $\left(3-\alpha\right)^{\text{th}},\alpha\in(0,1)$
order approximation to Caputo fractional (C-F) derivative using graded
mesh and standard central difference approximation for space derivatives,
in order to obtain the approximate solution of time fractional partial
differential equations (TFPDE). The proposed approximation for C-F
derivative tackles the singularity at origin effectively and is easily
applicable to diverse problems. The stability analysis and truncation
error bounds of the proposed scheme are discussed, along with this,
analyzed the required regularity of the solution. Few numerical examples
are presented to support the theory.
\end{abstract}

\section{Introduction\label{sec:Introduction HL1}}

\hspace{0.6 cm} The study of time fractional differential equations
is very intense in the past decade due to its applications in various
interdisciplinary areas, a detailed review has been presented in \cite{Ref Review}.
For the sake of presentation, in the following we study the constant
coefficient TFPDE, which is of the form:
\begin{equation}
_{0}^{C}D_{t}^{\alpha}y(x,t)=a\dfrac{\partial^{2}y(x,t)}{\partial x^{2}}+cy(x,t)+f(x,t),\quad\left(x,t\right)\in\left(0,X\right)\times(0,T],\,\alpha\in\left(0,1\right),a>0,c\le0\label{frde}
\end{equation}
\vspace{-0.5cm}
\begin{equation}
y(x,0)=\phi(x),\ x\in\left(0,X\right)\label{eq:frde ic}
\end{equation}
\vspace{-0.5cm}
\begin{equation}
y(0,t)=\psi_{0}(t),\;y(X,t)=\psi_{X}(t),\ t\in(0,T].\label{frde bc}
\end{equation}
Well adapted approximation of C-F derivative in the study of FDE is
the standard L1 scheme \cite{Ref Langlands,Jin ZHou} which is a $\left(2-\alpha\right)^{\text{th}}$
order approximation, the authors of \cite{JCAM Nagaraju and Harshita}
discussed an extension to standard L1 scheme which is a $\left(3-\alpha\right)^{\text{th}}$
order approximation. The authors in \cite{Ref Langlands,Jin ZHou,JCAM Nagaraju and Harshita},
while approximating the C-F derivative, did not consider the weak
singularity occurring at origin for the study of convergence of the
numerical scheme. Stynes et al., \cite{Ref Stynes SIAM} discussed
the numerical solution of TFPDE using standard L1 approximation on
graded mesh. This is one of the earliest articles to discuss the convergence
analysis of the method taking into consideration the initial weak
singularity. Very recently, Ren et al. in their article \cite{Ren appl math letters}
constructed a scheme based on the L1-type formula on graded mesh in
time and the direct discontinuous Galerkin in space directions for
solving TFPDE. Zheng et al. in \cite{Zheng appl math letters} presented
a scheme for solving the 2D multi-term time-fractional diffusion equation
with non-smooth solutions, where L1-type formula is derived on graded
mesh for approximating the C-F time derivative and Legendre spectral
approximation is used for the space derivatives.

In the following, we first state in brief the regularity requirement
for the solution of (\ref{frde})-(\ref{frde bc}), i.e., in section
2. In section 3, a numerical scheme is proposed based on the approximation
of the C-F derivative using second order non-uniform finite differences
and standard central differences for space derivatives. Stability
analysis and truncation error bounds are studied in section 4. Scrutinized
few examples for the applicability of the scheme in section 5.\vspace{-0.5cm}

\section{Regularity}

The series solution of TFPDE (\ref{frde})-(\ref{frde bc}) is well
discussed in \cite{Ref Luchko IVP} by using the variable separable
technique, the associated Strum-Liouville problem is\vspace{-0.5cm}
\[
\mathcal{L}p:=\left(-a\dfrac{d^{2}}{dx^{2}}-c\right)p=\lambda p;p(0)=0=p(X)
\]
Let $\lambda_{l}(>0)\text{ and }\theta_{l},l=1,2,\cdots$ be the eigen
values and normalised eigen functions respectively of this problem.

Based on the concepts of fractional sectorial operators, the domain
of $\mathcal{L^{\gamma}}$ is defined as\vspace{-0.25cm}
\[
D\left(\mathcal{L}^{\gamma}\right):=\left\{ f\in L_{2}\left(0,X\right):\sum_{l=1}^{\infty}\lambda_{l}^{2\gamma}\left|\left\langle f,\theta_{l}\right\rangle \right|^{2}<\infty\right\} .
\]
Let us define $\left\Vert f\right\Vert _{\mathcal{L}^{\gamma}}:=\left({\displaystyle \sum_{l=1}^{\infty}}\lambda_{l}^{2\gamma}\left|\left\langle f,\theta_{l}\right\rangle \right|^{2}\right)^{\frac{1}{2}}.$
Here $\left\langle .,.\right\rangle $ represents the standard scalar
product in $L_{2}\left(0,X\right).$

The series solution for (\ref{frde})-(\ref{frde bc}) with homogeneous
boundary conditions is given by
\[
y(x,t)=\sum_{l=1}^{\infty}\left[\left\langle \phi,\theta_{l}\right\rangle E_{\alpha,1}\left(-\lambda_{l}t^{\alpha}\right)+\int_{0}^{t}s^{\alpha-1}E_{\alpha,\alpha}\left(-\lambda_{l}t^{\alpha}\right)f_{l}(t-s)ds\right]\theta_{l}\left(x\right),
\]
where $f_{l}(s):=\left\langle f\left(\cdot,s\right),\theta_{l}\left(\cdot\right)\right\rangle $
and $E_{\alpha,\beta}$ is the Mittag-Leffler function defined by\vspace{-0.3cm}
\[
E_{\alpha,\beta}\left(s\right):={\displaystyle \sum_{k=0}^{\infty}}\dfrac{s^{k}}{\Gamma\left(\alpha k+\beta\right)}.
\]

In what follows, we state the theorem pertaining to the regularity
of the solution of (\ref{frde})-(\ref{frde bc}). The proof of this
theorem can be obtained on similar lines to the proof of Theorem 2.1
of \cite{Ref Stynes SIAM}.

\begin{The}\label{The: Regularity} Suppose $\phi\in D\left(\mathcal{L}^{7/2}\right)$,
$f\in D\left(\mathcal{L}^{5/2}\right),$ $f_{t},f_{tt}\in D\left(\mathcal{L}^{1/2}\right)$
and 
\[
\left\Vert f\left(.,t\right)\right\Vert _{\mathcal{L}^{5/2}}+\left\Vert f_{t}\left(.,t\right)\right\Vert _{\mathcal{L}^{1/2}}+\left\Vert f_{tt}\left(.,t\right)\right\Vert _{\mathcal{L}^{1/2}}+t^{\mu}\left\Vert f_{ttt}\left(.,t\right)\right\Vert _{\mathcal{L}^{1/2}}\leq K_{1},
\]
for all $t\in(0,T].$ Here, constant $K_{1}$ is not dependent on
$t$ and $\mu<1$ is an arbitrary constant. Then, the TFPDE with homogeneous
boundary conditions (\ref{frde})-(\ref{frde bc}) has a unique solution
$y(x,t)$ (satisfies the the differential equation and the initial
condition, point-wise), and $\exists\ K_{2}$ a constant $\ni$
\[
\begin{split}\left|\dfrac{\partial^{p}y}{\partial x^{p}}\left(x,t\right)\right| & \leq K_{2},\ p=0,1,2,3,4,\\
\left|\dfrac{\partial^{p}y}{\partial t^{p}}\left(x,t\right)\right| & \leq K_{2}\left(1+t^{\alpha-p}\right),\ p=0,1,2,3,
\end{split}
\ \forall\left(x,t\right)\in[0,X]\times(0,T].
\]
\end{The}

\section{Numerical scheme\label{sec:HL1}}

Suppose $M,N\in\mathbb{\mathbb{Z}^{+}}/\{0\}$ be the number of sub-intervals
in space and time direction respectively of the domain $\left(0,X\right)\times(0,T].$
The points in space are equidistant and in time direction are graded.
Let $\left(x_{m},t_{j}\right)$ be the discrete point in the domain,
we took $x_{0}=0,x_{m}=mh,x_{M}=X,t_{0}=0,t_{j}=\left(\frac{j}{N}\right)^{\beta}T,t_{N}=T.$
Here $\beta\in\mathbb{R^{+}}/\{0\}.$

First we take note on the approximation of C-F time derivative.\vspace{-0.2cm}

\subsection{Approximation of C-F derivative}

\hspace{0.6 cm}A higher order approximation for C-F derivative to
the function $u(t)$ is obtained using a modification to the standard
L1 scheme. From the definition of C-F derivative to $u(t)$ at $t=t_{j},$
we've\vspace{-0.2cm}
\begin{equation}
\begin{split}_{0}^{C}D_{t}^{\alpha}u(t_{j}) & =\frac{1}{\Gamma(1-\alpha)}{\displaystyle \int_{0}^{t_{j}}}(t_{j}-\eta)^{-\alpha}u^{\prime}(\eta)d\eta=\sum_{k=0}^{j-1}{\displaystyle \int_{t_{k}}^{t_{k+1}}}\dfrac{(t_{j}-\eta)^{-\alpha}}{\Gamma(1-\alpha)}u^{\prime}(\eta)d\eta=:\sum_{k=0}^{j-1}S_{k,j}.\end{split}
\label{eq: HL1 Caputo discretization}
\end{equation}
Considering the second order nonuniform finite difference approximation
for $u^{\prime}(t_{k}),u^{\prime\prime}(t_{k})$ in the Taylor's expansion
for $u^{\prime}$ taken as\vspace{-0.3cm}
\[
u^{\prime}(\eta)=u^{\prime}(t_{k})+\frac{(\eta-t_{k})}{1!}u^{\prime\prime}(t_{k})+\frac{(\eta-t_{k})^{2}}{2!}u^{\prime\prime\prime}(t_{k})+\mathcal{O}\left((\eta-t_{k})^{3}\right),\quad t_{k}\leq\eta\leq t_{k+1},
\]
upon simplifying we get\vspace{-0.2cm}
\begin{align*}
S_{k,j}= & \frac{1}{\Gamma(1-\alpha)}{\displaystyle \int_{t_{k}}^{t_{k+1}}}(t_{j}-\eta)^{-\alpha}\Bigg(\dfrac{u\left(t_{k+1}\right)-u\left(t_{k-1}\right)}{\tau_{k+1}+\tau_{k}}-\dfrac{\tau_{k+1}-\tau_{k}}{2}u^{\prime\prime}(t_{k})\\
 & +2(\eta-t_{k})\left[\dfrac{u\left(t_{k+1}\right)}{\tau_{k+1}\left(\tau_{k+1}+\tau_{k}\right)}-\dfrac{u\left(t_{k}\right)}{\tau_{k+1}\tau_{k}}+\dfrac{u\left(t_{k-1}\right)}{\tau_{k}\left(\tau_{k+1}+\tau_{k}\right)}\right]\Bigg)d\eta+Tr_{k,j},
\end{align*}
\begin{equation}
=2\delta_{k,j}\left[\dfrac{u^{k+1}}{\tau_{k+1}\left(\tau_{k+1}+\tau_{k}\right)}-\dfrac{u^{k}}{\tau_{k+1}\tau_{k}}+\dfrac{u^{k-1}}{\tau_{k}\left(\tau_{k+1}+\tau_{k}\right)}\right]+\gamma_{k,j}\left[\dfrac{u^{k+1}-u^{k-1}}{\tau_{k+1}+\tau_{k}}\right]+Tr_{k,j},\label{eq:L1 tj}
\end{equation}
where $Tr_{k,j}$ is the truncation error, $u^{k}\approxeq u\left(t_{k}\right)$,
$\tau_{k}=t_{k}-t_{k-1}$. Observe that for $k=0$ in equation (\ref{eq:L1 tj})
we lack the information for $u^{-1},$ to avoid this scenario $S_{0,j}$
is approximated separately. Initially for approximating $S_{0,1}$,
to be of order $\left(3-\alpha\right)$ the interval $\left[0,t_{1}\right]$
is divided into $\overline{N}$ sub intervals such that $\overline{N}^{\alpha-2}\approx N^{\alpha-3}t_{1},$
with $0=\overline{t}_{0}<\overline{t}_{1}<\overline{t}_{2}<\cdots<\overline{t}_{\overline{N}}=t_{1},$
and using the graded standard L1 approximation

\[
\begin{split}S_{0,1} & ={\displaystyle \sum_{k=0}^{\overline{N}-1}{\displaystyle \int_{\overline{t}_{k}}^{\overline{t}_{k+1}}}}\dfrac{(t_{1}-\eta)^{-\alpha}}{\Gamma(1-\alpha)}u^{\prime}(\eta)d\eta=-\zeta_{0,\overline{N}}u^{(0)}+\sum_{k=1}^{\overline{N}-1}u^{(k)}\left(\zeta_{k-1,\overline{N}}-\zeta_{k,\overline{N}}\right)+\zeta_{\overline{N}-1,\overline{N}}u^{(\overline{N})}+Tr_{1},\end{split}
\]
where $u^{(k)}\approxeq u\left(\overline{t}_{k}\right)$ and $\zeta_{k,\overline{N}}=\frac{\left(\overline{t}_{\overline{N}}-\overline{t}_{k}\right)^{1-\alpha}-\left(\overline{t}_{\overline{N}}-\overline{t}_{k+1}\right)^{1-\alpha}}{\Gamma(2-\alpha)\left(\overline{t}_{k+1}-\overline{t}_{k}\right)}.$\vspace{0.1cm}

Similarly, $S_{0,j},j\ge2$ occurring in (\ref{eq: HL1 Caputo discretization})
is also approximated using the graded standard L1 scheme as above.
Now, 
\begin{align*}
_{0}^{C}D_{t}^{\alpha}u(t_{j}) & =S_{0,j}+\sum_{k=1}^{j-1}S_{k,j}=-\xi_{0,j}u^{(0)}+\sum_{k=1}^{\overline{N}-1}u^{(k)}\left(\xi_{k-1,j}-\xi_{k,j}\right)+\xi_{\overline{N}-1,j}u^{(\overline{N})}\\
 & +\sum_{k=1}^{j-1}\left(\gamma_{k,j}\left[\dfrac{u^{k+1}-u^{k-1}}{\tau_{k+1}+\tau_{k}}\right]+2\delta_{k,j}\left[\dfrac{u^{k+1}}{\tau_{k+1}\left(\tau_{k+1}+\tau_{k}\right)}-\dfrac{u^{k}}{\tau_{k+1}\tau_{k}}+\dfrac{u^{k-1}}{\tau_{k}\left(\tau_{k+1}+\tau_{k}\right)}\right]\right)+Tr_{j},
\end{align*}
where $\xi_{k,j}=\frac{\left(t_{j}-\overline{t}_{k}\right)^{1-\alpha}-\left(t_{j}-\overline{t}_{k+1}\right)^{1-\alpha}}{\Gamma(2-\alpha)\left(\overline{t}_{k+1}-\overline{t}_{k}\right)}.$

\subsection{Derivation of the method}

Consider the discrete equation of (\ref{frde}) at $\left(x_{m},t_{j}\right)$
\[
_{0}^{C}D_{t}^{\alpha}y(x_{m},t_{j})=a\dfrac{\partial^{2}y(x_{m},t_{j})}{\partial x^{2}}+cy(x_{m},t_{j})+f(x_{m},t_{j}),\ 1\leq m\leq M-1,1\leq j\leq N.
\]
Replacing the approximation of the C-F derivative as discussed in
the previous subsection, central difference approximation for the
space derivative and upon simplification we get at $t=t_{1}$
\begin{equation}
\begin{split} & -y\left(x_{m-1},t_{1}\right)\dfrac{a}{h^{2}}+y\left(x_{m},t_{1}\right)\left[\zeta_{\overline{N}-1,\overline{N}}+\dfrac{2a}{h^{2}}-c\right]-y\left(x_{m+1},t_{1}\right)\dfrac{a}{h^{2}}\\
= & \zeta_{0,\overline{N}}y\left(x_{m},\overline{t}_{0}\right)+\sum_{k=1}^{\overline{N}-1}y\left(x_{m},\overline{t}_{k}\right)\left(\zeta_{k,\overline{N}}-\zeta_{k-1,\overline{N}}\right)+f\left(x,t_{1}\right)+Tr_{1}+\mathcal{O}\left(h^{2}\right)
\end{split}
\label{eq: HL1 num meth t1}
\end{equation}
and for $t=t_{j},\ j\geq2$
\begin{eqnarray}
 &  & -y\left(x_{m-1},t_{j}\right)\dfrac{a}{h^{2}}+y\left(x_{m},t_{j}\right)\left[d_{j,j}+\dfrac{2a}{h^{2}}-c\right]-y\left(x_{m+1},t_{j}\right)\dfrac{a}{h^{2}}\label{eq: HL1 num meth}\\
 & = & \xi_{0,j}y\left(x_{m},\overline{t}_{0}\right)+\sum_{k=1}^{\overline{N}-1}y\left(x_{m},\overline{t}_{k}\right)\left(\xi_{k,j}-\xi_{k-1,j}\right)-\xi_{\overline{N}-1,j}y\left(x_{m},\overline{t}_{\overline{N}}\right)+\sum_{k=0}^{j-1}d_{k,j}y\left(x_{m},t_{k}\right)+f\left(x_{m},t_{j}\right)+\mathcal{O}\left(h^{2}\right)+Tr_{j}.\nonumber 
\end{eqnarray}

\section{Stability analysis}

We discuss the stability of the proposed scheme for solving TFPDE
using Von-Neumann stability analysis. To do so, considered $\delta_{m}^{(J)}=z_{m}^{(J)}-y_{m}^{(J)}\text{ and }\ \delta_{m}^{j}=z_{m}^{j}-y_{m}^{j}$
$\ 0\leq m\leq M,\ 0\leq J\leq\overline{N}\ 0\leq j\leq N$ respectively
the difference between the perturbed and approximate solutions (perturbation
is due to a small variation in the initial condition) of the proposed
scheme as given in (\ref{eq: HL1 num meth t1}) and (\ref{eq: HL1 num meth}).
Replacing $\delta_{m}^{(J)}=\mu^{(J)}e^{i\rho mh}$ and $\delta_{m}^{j}=\mu^{j}e^{i\rho mh}$
where $\rho$ is the spatial wave number, $\mu$ is the amplitude
and $i^{2}=-1,$ in (\ref{eq: HL1 num meth t1}) and (\ref{eq: HL1 num meth})
respectively, we get\vspace{-0.2cm}
\begin{align}
 & \mu^{(\overline{N})}\overline{D}_{\overline{N}}=\zeta_{0,\overline{N}}\mu^{(0)}+\sum_{k=1}^{\overline{N}-1}\mu^{(k)}\left(\zeta_{k,\overline{N}}-\zeta_{k-1,\overline{N}}\right),\nonumber \\
 & \mu^{j}D_{j}=\sum_{k=0}^{j-1}d_{k,j}\mu^{k}+\xi_{0,j}\mu^{(0)}+\sum_{k=1}^{\overline{N}-1}\mu^{(k)}\left(\xi_{k,j}-\xi_{k-1,j}\right)-\xi_{\overline{N}-1,j}\mu^{(\overline{N})}.\label{eq: mu-j defn}
\end{align}
Further, for the stability of the scheme at $t=\overline{t}_{\overline{N}}=t_{1},$
one needs to analyze the intermediate calculations leading to (\ref{eq: HL1 num meth t1})
for which we have 
\begin{equation}
\mu^{(J)}\overline{D}_{J}=\zeta_{0,j}\mu^{(0)}+\sum_{k=1}^{J-1}\mu^{(k)}\left(\zeta_{k,J}-\zeta_{k-1,J}\right).\label{eq: mu J}
\end{equation}
Where $D_{j}=d_{j,j}-c+\dfrac{a\sin^{2}\left(\rho h/2\right)}{h^{2}},\ 1\leq j\leq N,$
$\overline{D}_{J}=\zeta_{J-1,J}-c+\dfrac{a\sin^{2}\left(\rho h/2\right)}{h^{2}},\ 1\leq J\leq\overline{N}.$

\vspace{0.2cm}
\hspace{-0.5cm}Note that $\overline{D}_{J}\geq\zeta_{J-1,J}\ \text{and\ }D_{j}\geq d_{j,j},$
holds true.\begin{Lem}\label{stability bounds HL1} $\left(a\right)$
For every $J\geq1,$ $\text{\ensuremath{\zeta_{k-1,J}}\ensuremath{\ensuremath{\leq\zeta_{k,J}},}\ \ensuremath{\forall k<J.} }$
$\left(b\right)$ For every $j\geq2,$ ${\displaystyle \sum_{k=0}^{j-1}}d_{k,j}=d_{j,j}.$

\end{Lem}

\begin{The} \label{Stability Thm} For every $1\leq J\leq\overline{N}$
and for $2\leq j\leq N$ respectively we have
\begin{equation}
\left|\mu^{(J)}\right|\leq\left|\mu^{0}\right|\text{ and }\left|\mu^{j}\right|\leq\left|\mu^{0}\right|.\label{eq: HL1 stability-1}
\end{equation}
\end{The}\bproof Define $\overline{\zeta}_{J}=\dfrac{\zeta_{0,J}+{\displaystyle \sum_{k=1}^{J-1}}\left(\zeta_{k,J}-\zeta_{k-1,J}\right)\overline{\zeta_{k}}}{\overline{D}_{J}},$
then equation (\ref{eq: mu J}) takes the form
\begin{equation}
\mu^{(J)}=\mu^{0}\overline{\zeta}_{J},\ 1\leq J\leq\overline{N.}\label{eq:mu j1}
\end{equation}
It is easy to see that to prove first part of (\ref{eq: HL1 stability-1})
we need to show that $\overline{\zeta}_{J}\leq1,$ for which we use
the principal of mathematical induction. Note that $\overline{\zeta}_{1}=\dfrac{\zeta_{0,1}}{\overline{D}_{1}}\leq1.$
For induction hypothesis, let $\overline{\zeta}_{J}\leq1,\ \forall J\leq K-1.$
Substituting this in the definition of $\overline{\zeta}_{J}$ at
$J=K$ and using lemma \ref{stability bounds HL1} gives
\[
\begin{split}\overline{\zeta}_{K}= & \dfrac{\zeta_{0,K}+{\displaystyle \sum_{k=1}^{K-1}}\left(\zeta_{k,K}-\zeta_{k-1,K}\right)\overline{\zeta}_{k}}{\overline{D}_{K}}\leq\dfrac{\zeta_{0,K}+{\displaystyle \sum_{k=1}^{K-1}}\left(\zeta_{k,K}-\zeta_{k-1,K}\right)}{\overline{D}_{K}}=\dfrac{\zeta_{K-1,K}}{\overline{D}_{K}}\leq1.\end{split}
\]
This implies $\overline{\zeta}_{J}\leq1,\ \forall1\leq J\leq\overline{N}.$

Now we prove the second part of (\ref{eq: HL1 stability-1}). From
equation (\ref{eq: mu-j defn}) we have
\begin{equation}
\begin{split}\mu^{j}= & \dfrac{{\displaystyle \sum_{k=0}^{j-1}}d_{k,j}\mu^{k}+\xi_{0,j}\mu^{(0)}+{\displaystyle \sum_{k=1}^{\overline{N}-1}}\mu^{(k)}\left(\xi_{k,j}-\xi_{k-1,j}\right)-\xi_{\overline{N}-1,j}\mu^{(\overline{N})}}{D_{j}}\\
= & \dfrac{d_{1,j}+{\displaystyle \sum_{k=2}^{j-1}}d_{k,j}\overline{d}_{1,k}-\xi_{\overline{N}-1,j}}{D_{j}}\mu^{1}+\dfrac{d_{0,j}+{\displaystyle \sum_{k=2}^{j-1}}d_{k,j}\overline{d}_{0,k}+\xi_{\overline{N}-1,j}}{D_{j}}\mu^{0}=:\overline{d}_{1,j}\mu^{1}+\overline{d}_{0,j}\mu^{0}.
\end{split}
\label{eq: muj in terms of mu1 and mu0}
\end{equation}
From the result in equation (\ref{eq:mu j1}) at $J=\overline{N},$
the equation (\ref{eq: muj in terms of mu1 and mu0}) reduces to
\begin{equation}
\mu^{j}\leq\left(\overline{d}_{1,j}+\overline{d}_{0,j}\right)\mu^{0}.\label{eq: muj in terms of mu0}
\end{equation}
One can see that from (\ref{eq: muj in terms of mu0}) the proof of
second part of (\ref{eq: HL1 stability-1}) follows by showing $\overline{d}_{1,j}+\overline{d}_{0,j}\leq1,$
for which we use the principal of mathematical induction. Note that
for $j=2$
\[
\overline{d}_{1,2}+\overline{d}_{0,2}=\dfrac{d_{0,2}+\xi_{\overline{N}-1,2}}{D_{2}}+\dfrac{d_{1,2}-\xi_{\overline{N}-1,2}}{D_{2}}=\dfrac{d_{2,2}}{D_{2}}\leq1.
\]
For induction hypothesis, let $\overline{d}_{1,j}+\overline{d}_{0,j}\leq1,\ \forall3\leq j\leq K-1.$
Substituting this in the definition of $\overline{d}_{0,K},\ \overline{d}_{1,K}$,
with the help of lemma \ref{stability bounds HL1} yields\vspace{-0.6cm}
\[
\begin{split}\overline{d}_{0,K}+\overline{d}_{1,K}= & \dfrac{d_{0,K}+{\displaystyle \sum_{k=2}^{K-1}}d_{k,K}\overline{d}_{0,k}-\xi_{\overline{N}-1,K}}{D_{K}}+\dfrac{d_{1,K}+{\displaystyle \sum_{k=2}^{K-1}}d_{k,K}\overline{d}_{1,k}+\xi_{\overline{N}-1,K}}{D_{K}}\\
= & \dfrac{d_{0,K}+d_{1,K}+{\displaystyle \sum_{k=2}^{K-1}}d_{k,K}\left(\overline{d}_{0,k}+\overline{d}_{1,k}\right)}{D_{K}}\leq\dfrac{{\displaystyle \sum_{k=0}^{K-1}}d_{k,K}}{D_{K}}=\dfrac{d_{K,K}}{D_{K}}\leq1.
\end{split}
\]
This implies $\overline{d}_{1,j}+\overline{d}_{0,j}\leq1,\ \forall1\leq j\leq N.$
\eproof\vspace{-0.6cm}

\section{Truncation error bounds}

The truncation error in time direction with $N$ mesh points is given
by
\begin{equation}
\left|Tr\right|=\left|\sum_{k=1}^{N-1}Tr_{k}\right|\leq\left|Tr_{1}\right|+\left|\sum_{k=2}^{N-1}{\displaystyle \int_{t_{k}}^{t_{k+1}}}\dfrac{(t_{j}-x)^{-\alpha}}{\Gamma\left(1-\alpha\right)}y^{\prime\prime\prime}(t_{k})\left\{ \frac{(x-t_{k})^{2}}{2}-(x-t_{k})\dfrac{\tau_{k+1}-\tau_{k}}{3}-\dfrac{\tau_{k+1}\tau_{k}}{6}\right\} dx\right|.\label{eq: HL1 Tr}
\end{equation}
Simplifying the above equation and from Theorem \ref{The: Regularity}
we have
\begin{equation}
\begin{split}\left|Tr\right|\leq & \sum_{k=2}^{N-1}\dfrac{\left(k^{\beta}\right)^{\alpha-3}}{36\Gamma\left(4-\alpha\right)}\Bigg|\left(N^{\beta}-k^{\beta}\right)^{1-\alpha}\biggl\{\left(\alpha-2\right)\left(\alpha-3\right)\left(k+1\right)^{\beta}\left(k-1\right)^{\beta}-\alpha\left(1-\alpha\right)k^{2\beta}\\
 & +\alpha\left(3-\alpha\right)k^{\beta}\left(\left(k+1\right)^{\beta}+\left(k-1\right)^{\beta}\right)-2N^{\beta}\left(\left(3-\alpha\right)\left(\left(k+1\right)^{\beta}+\left(k-1\right)^{\beta}\right)+2k^{\beta}-3N^{\beta}\right)\biggr\}\\
 & -\left(N^{\beta}-\left(k+1\right)^{\beta}\right)^{1-\alpha}\biggl\{\left(\alpha-2\right)\left(\alpha-3\right)k^{\beta}\left(k-1\right)^{\beta}+\alpha\left(3-\alpha\right)\left(k+1\right)^{\beta}\left(k^{\beta}+\left(k-1\right)^{\beta}\right)\\
 & -\alpha\left(1-\alpha\right)\left(k+1\right)^{2\beta}-2N^{\beta}\left(\left(3-\alpha\right)\left(k^{\beta}+\left(k-1\right)^{\beta}\right)+2\left(k+1\right)^{\beta}-3N^{\beta}\right)\biggr\}\Bigg|+\min\left\{ \overline{N}^{\alpha-2},\overline{N}^{-\beta\alpha}\right\} \\
\leq & \min\left\{ N^{\alpha-3},N^{-\beta\alpha}\right\} .
\end{split}
\label{eq: Trunc eror HL1}
\end{equation}
As a consequence of equation (\ref{eq: Trunc eror HL1}) along with
the theorem \ref{Stability Thm} one can have \begin{The}\label{Convergence Thm HL1}
The solution of numerical scheme $y_{i}^{j},$ satisfies
\begin{equation}
\left|y_{i}^{j}-y\left(x_{i},t_{j}\right)\right|\leq C\left(h^{2}+N^{\min\left\{ \alpha-3,-\beta\alpha\right\} }\right)\label{eq: HL1 convergence thm}
\end{equation}
 \end{The}

\section{Numerical illustrations\label{sec:Numerics HL1}}

In this section we present three diverse examples: First a fractional
delay differential equation with non-smooth solution, then a time
fractional diffusion equation and nonlinear TFPDE, to understand the
applicability of the proposed method. All the examples exhibit weak
initial singularity. A comparison between the proposed scheme and
the L1 scheme is also shown. In all the following results the optimal
value for $\beta$ is considered.

\begin{figure}[H]
\begin{minipage}[c][1\totalheight][t]{0.45\textwidth}%
\captionof{table}{$L_\infty$ errors of example \ref{eg 2 fdde}}
\begin{flushleft}
{\scriptsize{}}%
\begin{tabular}{lclcl}
\hline 
\multirow{2}{*}{{\scriptsize{}$N$}} & \multicolumn{2}{c}{{\scriptsize{}Uniform mesh}} & \multicolumn{2}{c}{{\scriptsize{}Graded mesh with $\beta=5$}}\tabularnewline
 & {\scriptsize{}$L_{\infty}$ error} & {\scriptsize{}EOC} & {\scriptsize{}$L_{\infty}$ error} & {\scriptsize{}EOC}\tabularnewline
\hline 
{\scriptsize{}$2^{6}$} & {\scriptsize{}$8.90\times10^{-3}$} &  & {\scriptsize{}$3.61\times10^{-3}$} & \tabularnewline
{\scriptsize{}$2^{7}$} & {\scriptsize{}$7.31\times10^{-3}$} & {\scriptsize{}$0.284$} & {\scriptsize{}$6.39\times10^{-4}$} & {\scriptsize{}$2.500$}\tabularnewline
{\scriptsize{}$2^{8}$} & {\scriptsize{}$5.63\times10^{-3}$} & {\scriptsize{}$0.377$} & {\scriptsize{}$1.13\times10^{-4}$} & {\scriptsize{}$2.500$}\tabularnewline
{\scriptsize{}$2^{9}$} & {\scriptsize{}$4.21\times10^{-3}$} & {\scriptsize{}$0.420$} & {\scriptsize{}$1.99\times10^{-5}$} & {\scriptsize{}$2.500$}\tabularnewline
{\scriptsize{}$2^{10}$} & {\scriptsize{}$3.08\times10^{-3}$} & {\scriptsize{}$0.451$} & {\scriptsize{}$3.53\times10^{-6}$} & {\scriptsize{}$2.500$}\tabularnewline
\hline 
\end{tabular}{\scriptsize\par}
\par\end{flushleft}%
\end{minipage}\hfill{}%
\begin{minipage}[c][1\totalheight][t]{0.45\textwidth}%
\includegraphics[scale=0.4]{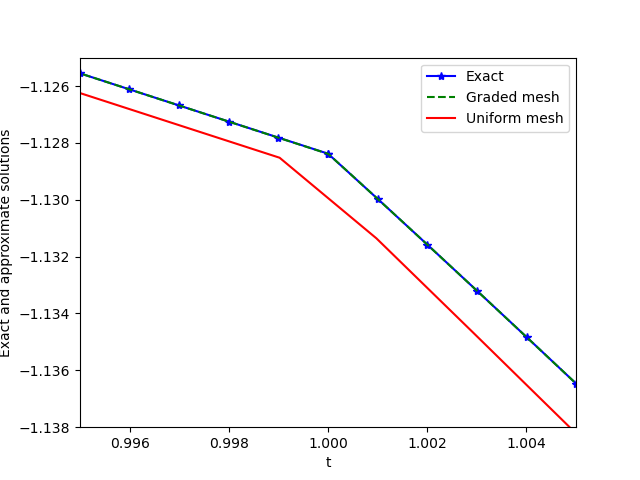}\vspace*{-0.25cm}

\caption{\label{fig: fdde}Plot for the results in Ex. \ref{eg 2 fdde}}
\end{minipage}
\end{figure}

\begin{Exp}\label{eg 2 fdde}Consider the following fractional delay
differential equation with nonsmooth solution
\[
_{0}^{C}D_{t}^{\frac{1}{2}}y(t)=y(t-1)-t,\quad t\in\left[0,2\right]
\]
\[
y(t)=t,\quad t\in[-1,0]
\]
whose analytical solution is given in detail in \cite{Ref Morgado}\end{Exp}$L_{\infty}$
errors obtained by HL1 scheme for the of the above example with non-smooth
solution are given in table 1 which shows that even with few mesh
points the error obtained using graded mesh is much better than uniform
mesh. The exact and approximate solutions at $t=1$ (where the solution
is not smooth) are plotted in the figure \ref{fig: fdde}, it can
be seen that the HL1 scheme gives a good resolution.\vspace{-0.4cm}

\begin{figure}[H]
\qquad{}\qquad{}\includegraphics[width=15cm,height=5cm]{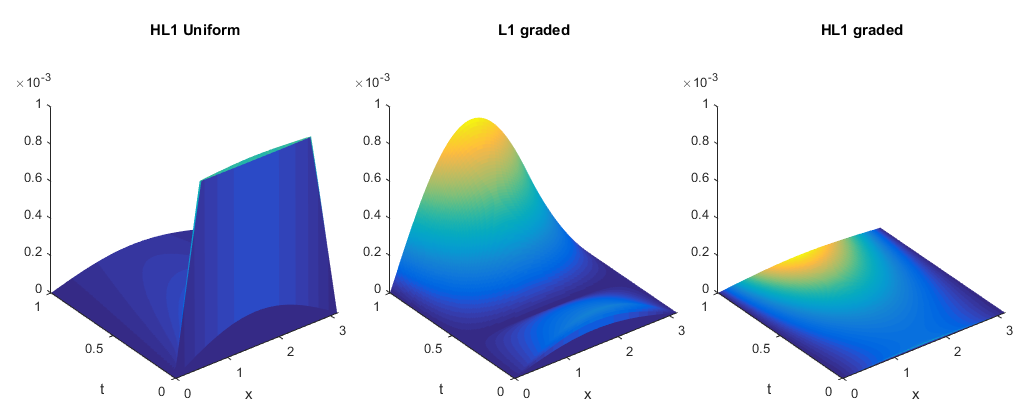}\caption{\label{fig: fpde abs err}Absolute point-wise error plots for example
\ref{eg fpde}}
\end{figure}

\begin{Exp} \label{eg fpde}Consider the fractional diffusion equation
of the form (\ref{frde})-(\ref{frde bc}) with $a=1,c=0,X=\pi,T=1.$
whose exact solution is given by $y(x,t)=\left(t^{3}+t^{\alpha}\right)\sin(x).$\end{Exp}
The tables \ref{tab: fpde uniform-HL1}, \ref{tab:fpde graded-HL1}
display the maximum absolute errors (MAE) for this example . It is
clear from these tables that the numerical results match the theoretical
estimates for graded meshes. A comparison between L1 and HL1 in graded
meshes shows that HL1 gives better results. This can also be understood
through figure \ref{fig: fpde abs err}( for $M=256=N$ and $\alpha=0.6$).
One can observe from this figure that high absolute error occurring
at $t=0$ diminishes significantly using the proposed scheme. \vspace{-0.5cm}
\begin{table}[H]
\centering{}{\small{}\caption{\label{tab: fpde uniform-HL1}Maximum absolute errors for example
\ref{eg fpde} using HL1 scheme on uniform mesh}
}{\scriptsize{}}%
\begin{tabular}{lcccccc}
\hline 
\multirow{2}{*}{{\scriptsize{}$M=N$}} & \multicolumn{2}{c}{{\scriptsize{}$\alpha=0.4$}} & \multicolumn{2}{c}{{\scriptsize{}$\alpha=0.6$}} & \multicolumn{2}{c}{{\scriptsize{}$\alpha=0.8$}}\tabularnewline
 & {\scriptsize{}MAE} & {\scriptsize{}EOC} & {\scriptsize{}MAE} & {\scriptsize{}EOC} & {\scriptsize{}MAE} & {\scriptsize{}EOC}\tabularnewline
\hline 
\multirow{1}{*}{{\scriptsize{}$2^{5}$}} & {\scriptsize{}$1.04\times10^{-2}$} &  & {\scriptsize{}$5.45\times10^{-3}$} &  & {\scriptsize{}$2.42\times10^{-3}$} & \tabularnewline
\multirow{1}{*}{{\scriptsize{}$2^{6}$}} & {\scriptsize{}$8.69\times10^{-3}$} & {\scriptsize{}$0.26$} & {\scriptsize{}$4.14\times10^{-3}$} & {\scriptsize{}$0.40$} & {\scriptsize{}$1.58\times10^{-3}$} & {\scriptsize{}$0.61$}\tabularnewline
\multirow{1}{*}{{\scriptsize{}$2^{7}$}} & {\scriptsize{}$7.21\times10^{-3}$} & {\scriptsize{}$0.27$} & {\scriptsize{}$3.00\times10^{-3}$} & {\scriptsize{}$0.46$} & {\scriptsize{}$9.99\times10^{-4}$} & {\scriptsize{}$0.66$}\tabularnewline
\multirow{1}{*}{{\scriptsize{}$2^{8}$}} & {\scriptsize{}$5.79\times10^{-3}$} & {\scriptsize{}$0.31$} & {\scriptsize{}$2.11\times10^{-3}$} & {\scriptsize{}$0.50$} & {\scriptsize{}$6.08\times10^{-4}$} & {\scriptsize{}$0.72$}\tabularnewline
\multirow{1}{*}{{\scriptsize{}$2^{9}$}} & {\scriptsize{}$4.59\times10^{-3}$} & {\scriptsize{}$0.34$} & {\scriptsize{}$1.45\times10^{-3}$} & {\scriptsize{}$0.54$} & {\scriptsize{}$3.61\times10^{-4}$} & {\scriptsize{}$0.75$}\tabularnewline
\hline 
\end{tabular}
\end{table}
\vspace{-0.5cm}
\begin{table}[H]
\centering{}{\small{}\caption{\label{tab:fpde graded-HL1}Maximum absolute errors for example \ref{eg fpde}
using}
}{\scriptsize{}}%
\begin{tabular}{lcccccc|}
\hline 
\multirow{3}{*}{{\scriptsize{}$M=N$}} & \multicolumn{6}{c}{{\scriptsize{}HL1 scheme with $\beta=(3-\alpha)/\alpha$}}\tabularnewline
\cline{2-7} \cline{3-7} \cline{4-7} \cline{5-7} \cline{6-7} \cline{7-7} 
 & \multicolumn{2}{c}{{\scriptsize{}$\alpha=0.4$}} & \multicolumn{2}{c}{{\scriptsize{}$\alpha=0.6$}} & \multicolumn{2}{c|}{{\scriptsize{}$\alpha=0.8$}}\tabularnewline
 & {\scriptsize{}MAE} & {\scriptsize{}EOC} & {\scriptsize{}MAE} & {\scriptsize{}EOC} & {\scriptsize{}MAE} & {\scriptsize{}EOC}\tabularnewline
\hline 
\multirow{1}{*}{{\scriptsize{}$2^{6}$}} & {\scriptsize{}$1.35\times10^{-3}$} &  & {\scriptsize{}$6.43\times10^{-4}$} &  & {\scriptsize{}$7.00\times10^{-4}$} & \tabularnewline
\multirow{1}{*}{{\scriptsize{}$2^{7}$}} & {\scriptsize{}$2.24\times10^{-4}$} & {\scriptsize{}$2.59$} & {\scriptsize{}$1.36\times10^{-4}$} & {\scriptsize{}$2.24$} & {\scriptsize{}$1.60\times10^{-4}$} & {\scriptsize{}$2.13$}\tabularnewline
\multirow{1}{*}{{\scriptsize{}$2^{8}$}} & {\scriptsize{}$3.69\times10^{-5}$} & {\scriptsize{}$2.60$} & {\scriptsize{}$2.86\times10^{-5}$} & {\scriptsize{}$2.24$} & {\scriptsize{}$3.62\times10^{-5}$} & {\scriptsize{}$2.14$}\tabularnewline
\multirow{1}{*}{{\scriptsize{}$2^{9}$}} & {\scriptsize{}$6.09\times10^{-6}$} & {\scriptsize{}$2.60$} & {\scriptsize{}$6.08\times10^{-6}$} & {\scriptsize{}$2.25$} & {\scriptsize{}$8.20\times10^{-6}$} & {\scriptsize{}$2.14$}\tabularnewline
\multirow{1}{*}{{\scriptsize{}$2^{10}$}} & {\scriptsize{}$1.11\times10^{-6}$} & {\scriptsize{}$2.60$} & {\scriptsize{}$1.31\times10^{-6}$} & {\scriptsize{}$2.25$} & {\scriptsize{}$1.86\times10^{-6}$} & {\scriptsize{}$2.14$}\tabularnewline
\hline 
\end{tabular}{\scriptsize{}}%
\begin{tabular}{|cccccc}
\hline 
\multicolumn{6}{c}{{\scriptsize{}L1 scheme with $\beta=(2-\alpha)/\alpha$}}\tabularnewline
\hline 
\multicolumn{2}{|c}{{\scriptsize{}$\alpha=0.4$}} & \multicolumn{2}{c}{{\scriptsize{}$\alpha=0.6$}} & \multicolumn{2}{c}{{\scriptsize{}$\alpha=0.8$}}\tabularnewline
{\scriptsize{}MAE} & {\scriptsize{}EOC} & {\scriptsize{}MAE} & {\scriptsize{}EOC} & {\scriptsize{}MAE} & {\scriptsize{}EOC}\tabularnewline
\hline 
{\scriptsize{}$4.14\times10^{-3}$} &  & {\scriptsize{}$5.12\times10^{-3}$} &  & {\scriptsize{}$7.98\times10^{-3}$} & \tabularnewline
{\scriptsize{}$1.45\times10^{-3}$} & {\scriptsize{}$1.51$} & {\scriptsize{}$1.98\times10^{-3}$} & {\scriptsize{}$1.37$} & {\scriptsize{}$3.47\times10^{-3}$} & {\scriptsize{}$1.20$}\tabularnewline
{\scriptsize{}$5.01\times10^{-4}$} & {\scriptsize{}$1.54$} & {\scriptsize{}$7.60\times10^{-4}$} & {\scriptsize{}$1.38$} & {\scriptsize{}$1.51\times10^{-3}$} & {\scriptsize{}$1.20$}\tabularnewline
{\scriptsize{}$1.70\times10^{-4}$} & {\scriptsize{}$1.55$} & {\scriptsize{}$2.90\times10^{-4}$} & {\scriptsize{}$1.39$} & {\scriptsize{}$6.56\times10^{-4}$} & {\scriptsize{}$1.20$}\tabularnewline
{\scriptsize{}$5.75\times10^{-5}$} & {\scriptsize{}$1.56$} & {\scriptsize{}$1.11\times10^{-4}$} & {\scriptsize{}$1.39$} & {\scriptsize{}$2.85\times10^{-4}$} & {\scriptsize{}$1.20$}\tabularnewline
\hline 
\end{tabular}
\end{table}
\vspace{-0.5cm}

\begin{Exp} \label{eg 5 nlpde HL1}Consider the following nonlinear
TFPDE \cite{Ref Baranwal CNSNS 2012}
\[
_{0}^{C}D_{t}^{\alpha}y(x,t)=\dfrac{\partial^{2}y(x,t)}{\partial x^{2}}-y(x,t)\left(1-y(x,t)\right)+f(x,t),x\in(0,1),t>0,
\]
with side conditions
\[
y(x,0)=1,x\in[0,1];\quad y(0,t)=1,\;y(1,t)=1+\dfrac{t^{\alpha}}{\Gamma\left(1+\alpha\right)}\sin(1),t\ge0.
\]
\end{Exp}The source term is evaluated by taking $y(x,t)=1+\dfrac{t^{\alpha}}{\Gamma\left(1+\alpha\right)}\sin(x),$
as an exact solution. Here the nonlinear example is converted into
linearized system of equations using Newton's quasi-linearization
method. Table \ref{tab: nonlinear pde graded HL1} displays the maximum
absolute errors for this example at $\alpha=0.4,0.6,0.8$ and it can
be seen that the expected order of convergence (EOC) is achieved.
This example was discussed to show the applicability of HL1 scheme
to nonlinear problems. \vspace{-0.5cm}
\begin{table}[H]
\centering{}{\small{}\caption{\label{tab: nonlinear pde graded HL1}Maximum absolute errors for
example \ref{eg 5 nlpde HL1} using graded mesh with $\beta=(3-\alpha)/\alpha$}
}{\scriptsize{}}%
\begin{tabular}{lcccccc}
\hline 
\multirow{2}{*}{{\scriptsize{}$M=N$}} & \multicolumn{2}{c}{{\scriptsize{}$\alpha=0.4$}} & \multicolumn{2}{c}{{\scriptsize{}$\alpha=0.6$}} & \multicolumn{2}{c}{{\scriptsize{}$\alpha=0.8$}}\tabularnewline
 & {\scriptsize{}MAE} & {\scriptsize{}EOC} & {\scriptsize{}MAE} & {\scriptsize{}EOC} & {\scriptsize{}MAE} & {\scriptsize{}EOC}\tabularnewline
\hline 
\multirow{1}{*}{{\scriptsize{}$2^{4}$}} & {\scriptsize{}$5.28\times10^{-2}$} &  & {\scriptsize{}$1.05\times10^{-2}$} &  & {\scriptsize{}$3.47\times10^{-3}$} & \tabularnewline
\multirow{1}{*}{{\scriptsize{}$2^{5}$}} & {\scriptsize{}$9.18\times10^{-3}$} & {\scriptsize{}$2.524$} & {\scriptsize{}$2.11\times10^{-3}$} & {\scriptsize{}$2.254$} & {\scriptsize{}$9.06\times10^{-4}$} & {\scriptsize{}$1.940$}\tabularnewline
\multirow{1}{*}{{\scriptsize{}$2^{6}$}} & {\scriptsize{}$1.53\times10^{-3}$} & {\scriptsize{}$2.588$} & {\scriptsize{}$4.87\times10^{-4}$} & {\scriptsize{}$2.372$} & {\scriptsize{}$2.14\times10^{-4}$} & {\scriptsize{}$2.081$}\tabularnewline
\multirow{1}{*}{{\scriptsize{}$2^{7}$}} & {\scriptsize{}$2.52\times10^{-4}$} & {\scriptsize{}$2.598$} & {\scriptsize{}$7.73\times10^{-5}$} & {\scriptsize{}$2.395$} & {\scriptsize{}$4.86\times10^{-5}$} & {\scriptsize{}$2.139$}\tabularnewline
\multirow{1}{*}{{\scriptsize{}$2^{8}$}} & {\scriptsize{}$4.16\times10^{-5}$} & {\scriptsize{}$2.600$} & {\scriptsize{}$1.46\times10^{-5}$} & {\scriptsize{}$2.399$} & {\scriptsize{}$1.08\times10^{-5}$} & {\scriptsize{}$2.170$}\tabularnewline
\hline 
\end{tabular}
\end{table}

\section{Conclusion\label{sec:Conclusion}}

\hspace{0.6 cm}In this article, we presented the HL1 scheme on graded
mesh ($\beta$ the mesh ratio) by taking into consideration the initial
singularity arising in the time fractional derivative. Stability analysis
and truncation error bounds for the proposed scheme are discussed.
The scheme using graded mesh has the order of accuracy to be $\min\left\{ \beta\alpha,3-\alpha\right\} .$
It is evident from the numerical examples that the graded mesh scheme
resolves the singularity with high resolution by attaining the desired
order of accuracy.

\end{document}